\documentclass[12pt]{amsart}
\usepackage{amsmath, amsthm, amscd, amssymb, amsfonts, latexsym}
\usepackage{fullpage}
\usepackage{bm,mathdots}
\usepackage{dsfont}
\usepackage[all]{xy}
\usepackage{tikz}
\usepackage{fancybox}
\usepackage{ytableau}

\newcommand{\kommentar}[1]{}

\newcommand{\F}{\mathbb F}

\newcommand{\Z}{\mathbb Z}

\newcommand{\R}{\mathbb R}

\newcommand{\C}{\mathbb C}

\DeclareMathOperator{\re}{Re}
\DeclareMathOperator{\im}{Im}

\renewcommand{\pmod}[1]{\,(\mathrm{mod}\,#1)}

\newtheorem{lem}{Lemma}[section]
\newtheorem{prop}[lem]{Proposition}
\newtheorem{thm}[lem]{Theorem}

\newtheorem{conj}[lem]{Conjecture}
\theoremstyle{definition}

\newtheorem{rem}[lem]{Remark}

\author{Vivian Kuperberg}
\author{Matilde Lal\'in}
\address{Vivian Kuperberg: ETH Z\"urich - Departement Mathematik,
R\"amistrasse 101,
8092 Z\"urich,
Switzerland} \email{vivian.kuperberg@math.ethz.ch}

\address{Matilde Lal\'in:  D\'epartement de math\'ematiques et de statistique, Universit\'e de Montr\'eal, CP 6128, succ. Centre-ville, Montreal, QC H3C 3J7, Canada}\email{matilde.lalin@umontreal.ca}

\thanks{This work is supported by NSF GRFP grant DGE-1656518, the NSF Mathematical Sciences Research Program through the grant DMS-2202128,  the Natural Sciences and Engineering Research Council of Canada, RGPIN-2022-03651, and the Fonds de recherche du Qu\'ebec - Nature et technologies, Projet de recherche en \'equipe 300951 and 345672. V.K. would also like to thank the Mittag-Leffler Institute for its hospitality.}


\title{Arithmetic constants for symplectic variances of the divisor function}

\begin{document}

\begin{abstract}
In \cite{KuperbergLalin2}, the authors stated some conjectures on the variance of certain sums of the divisor function $d_k(n)$ over number fields, which were inspired by analogous results over function fields proven in \cite{KuperbergLalin}. These problems are related to certain symplectic matrix integrals. While the function field results can be directly related to the random matrix integrals, the connection between the random matrix integrals and the number field results is less direct and involves arithmetic factors. The goal of this article is to give heuristic arguments for the formulas of these arithmetic factors.
\end{abstract}

\maketitle

\section{Introduction}

The goal of this article is to investigate the arithmetic constants appearing in certain conjectures made in \cite{KuperbergLalin2} involving the divisor function and integrals over the ensemble of unitary symplectic matrices. For a positive integer $n$, let $d_k(n)$ denote the number of ways of writing $n$ as an ordered product of $k$ positive integers. In \cite{KR3}, Keating, Rodgers, Roditty-Gershon, and Rudnick study the distribution of the sum
\begin{equation} \label{eq:Snf}
\mathcal{S}_{d_k;X;Q}(A):= \sum_{\substack{m\leq X\\m\equiv A \pmod{Q}}}d_k(m),\end{equation}
where $Q\in \Z$ is prime, and $A\in \Z$. After proving the analogous result for a function field model of the same problem, they conjecture in \cite[Conjecture 3.3]{KR3} that when $Q^{1+\varepsilon}<X<Q^{k-\varepsilon}$ and as $X\rightarrow \infty$,
 the variance of $\mathcal{S}_{d_k;X;Q}$, defined by
\[\mathrm{Var}(\mathcal{S}_{d_k;X;Q}):=\frac{1}{\Phi(Q)} \sum_{\substack{A\pmod{Q}\\\mathrm{gcd}(A,Q)=1}} \left|\mathcal{S}_{d_k;X;Q}(A)-\left \langle \mathcal{S}_{d_k;X;Q}\right\rangle  \right|^2,\]
where
\[\left \langle \mathcal{S}_{d_k;X;Q}\right\rangle:= \frac{1}{\Phi(Q)} \sum_{\substack{A\pmod{Q}\\\mathrm{gcd}(A,Q)=1}}\mathcal{S}_{d_k;X;Q}(A),\]
is given by
\begin{equation} \label{eq:KR3-Conj5.3}
\mathrm{Var}(\mathcal{S}_{d_k;X;Q}) \sim \frac{X}{Q} a_k \gamma_k\left(\frac{\log X}{\log Q}\right)(\log Q)^{k^2-1}.
\end{equation}
Here
\[a_k:=\prod_p \left( \Big(1-\frac{1}{p}\Big)^{(k-1)^2}\sum_{j=0}^{k-1}
\binom{k-1}{j}^2\frac{1}{p^j}\right)\]
 and $\gamma_k(c)$ is a piecewise polynomial of degree $k^2-1$ defined as the leading coefficient (respect to the variable $N$) of $I_{k}(cN;N)$, where
\[
I_{k}(m;N):=\int_{\mathrm{U}(N)} \Big|\sum_{\substack{j_1+\cdots+j_k=m\\0\leq j_1,\dots,j_k \leq N}}\mathrm{Sc}_{j_1}(U)\cdots \mathrm{Sc}_{j_k}(U)\Big|^2 \mathrm{d}U.
\]
The integral ranges over complex unitary matrices of dimension $N$, and the $\mathrm{Sc}_j(U)$ are the secular coefficients, defined for a $N \times N$ matrix $U$ via
\[\det (I+Ux)=\sum_{j=0}^N \mathrm{Sc}_j(U)x^j.\]

More precisely, letting $c=\frac{m}{N} \in [0,k]$, Keating et al. further prove in \cite[Theorem 1.5]{KR3} that
 \[\int_{\mathrm{U}(N)}\Big| \sum_{\substack{j_1+\cdots+j_k=m\\0\leq j_1,\dots,j_k\leq N}}\mathrm{Sc}_{j_1}(U)\cdots \mathrm{Sc}_{j_k}(U)\Big|^2 \mathrm{d}U=\gamma_k(c)N^{k^2-1}+O_k(N^{k^2-2}).\]

The conjecture given by equation \eqref{eq:KR3-Conj5.3} follows from \cite[Theorem 3.1]{KR3}, which is a function field version of the same problem. Let $\F_q[T]$ be the ring of polynomials over the function field $\F_q$ of $q$ elements, where $q$ is an odd prime power. For $f\in \F_q[T]$, we define  $d_k(f)$ as the number of ways of writing $f$ as an ordered product of $k$ monic polynomials. For $Q\in \F_q[T]$ square-free and $A\in \F_q[T]$ coprime to $Q$,  Keating et al consider the sum
\begin{equation}\label{eq:Sff}
\mathcal{S}_{d_k;n;Q}(A):= \sum_{\substack{f \text{monic}\\ \deg(f)=n\\f\equiv A \pmod{Q}}}d_k(f),\end{equation}
and prove that for $n \leq k (\deg (Q)-1)$ and as $q \rightarrow \infty$, the variance of the sum above is given by
\begin{equation} \label{eq:KR3-Thm3.1}
\mathrm{Var}(\mathcal{S}_{d_k;n;Q}) \sim \frac{q^n}{q^{\deg(Q)}} I_{k}(n;\deg(Q)-1)\sim \frac{q^n}{q^{\deg(Q)}}
\gamma_k\left(\frac{n}{\deg(Q)}\right)(\deg(Q))^{k^2-1}.
\end{equation}
To understand the relationship between the result \eqref{eq:KR3-Thm3.1} and the conjectural statement \eqref{eq:KR3-Conj5.3}, we observe that $q^n$ represents the size of the sum over $f$ in \eqref{eq:Sff}, and it therefore corresponds to $X$, the size of the sum over $m$ in \eqref{eq:Snf}. Similarly, $q^{\deg(Q)}$ represents the size of $Q$, and $n$ represents $\log X$, while $\deg(Q)$ represents $\log(Q)$. (The statement of the conjectural formula \eqref{eq:KR3-Conj5.3} is done with the assumption $Q$ prime for simplicity.)

The key (and only!) difference between the function field result that Keating et al. prove in \cite[Theorem 3.1]{KR3} and \cite[Conjecture 3.3]{KR3} is the arithmetic factor $a_k$. This term was computed in a conjecture of K\"osters \cite{Kosters} relating an integral of shifted zeta values to an analogous random matrix integral. The integral of the shifted zeta values is very close to the quantity that is conjecturally computed in the integer case in \cite{KR3}, whereas the random matrix integral is very close to the same quantity in the function field case, so the relation between the two is precisely the arithmetic factor $a_k$. The same arithmetic factor also appears in work of Lester \cite{Lester} on the variances of divisor functions on middle length intervals.

Keating et al. also propose a similar conjecture for the variance over short intervals in \cite[Conjecture 1.1]{KR3}, based on a function field version \cite[Theorem 1.2]{KR3}.

In \cite{KuperbergLalin2} the authors of this article formulate some conjectures over number fields based on their previous work in \cite{KuperbergLalin} over function fields. A goal of both articles \cite{KuperbergLalin,KuperbergLalin2} is to study similar problems to those considered in \cite{KR3}, but leading to integrals involving
the group of unitary symplectic matrices.

\begin{conj}\cite[Conjecture 1.10]{KuperbergLalin2} \label{conj:symp-square-fund-discs} For a positive fundamental discriminant $r$, define
\[\mathcal{T}^S_{d_k;x}(r)=\sum_{\substack{n\leq x  \\(n,r) = 1}}d_k(n)\chi_r(n),\]
where $\chi_r$ is the primitive quadratic character mod $r$.
Let $y > 0$. Define the variance of $\mathcal T^S_{d_k;x}(r)$ in $(y,2y]$ by
\begin{align*}
\mathrm{Var}_{r\in (y,2y]}\left(\mathcal{T}^S_{d_k;x}\right) &:= \mathbb E^*_{y < r \le 2y} \Big(\mathcal T^S_{d_k;x}(r) \Big)^2,
\end{align*}
where the expectation $\mathbb E^*$ is taken over positive fundamental discriminants.

Then, for $x^{1/k+\varepsilon}\leq y$, we have,
\begin{equation}\label{eq:conj-fund-discs}
\mathrm{Var}_{r\in (y, 2y]}\left(\mathcal{T}^S_{d_k;x}\right)
 \sim a_k^S(\mathcal T) x\gamma_{d_k,2}^S\left(\frac{\log x}{\log y}\right) (\log y)^{2k^2+k-2},
\end{equation}
where $a_k^S(\mathcal T)$ is an arithmetic constant given by
\begin{equation} \label{eq:akS}
a_k^S(\mathcal T) = 2 \prod_p \left(1-\frac 1p\right)^{k(2k+1)} \left(\frac 1{p+1}\left(1 + \frac p2\left(\left(1 + \frac 1{\sqrt p}\right)^{-2k} + \left(1 - \frac 1{\sqrt p}\right)^{-2k}\right)\right)\right),
\end{equation}
and $\gamma_{d_k,2}^S(c)$ is a piecewise polynomial of degree $2k^2+k-2$ defined as the leading coefficient (respect to the variable $N$) of $I_{d_k,2}^S(c2N;N)$, where
\[
I_{d_k,2}^S(m;N):=\int_{\mathrm{Sp}(2N)} \Big|\sum_{\substack{j_1+\cdots+j_k=m\\0\leq j_1,\dots,j_k \leq 2N}}\mathrm{Sc}_{j_1}(U)\cdots \mathrm{Sc}_{j_k}(U)\Big|^2 \mathrm{d}U.
\]
\end{conj}

As in the work of \cite{KR3}, Conjecture \ref{conj:symp-square-fund-discs} is based on a function field version that contains analogous elements to all the factors on the right-hand side of \eqref{eq:conj-fund-discs}, with the exception of the arithmetic factor $a_k^S(\mathcal T)$. A goal of the present article is to give a heuristic argument of why we expect $a_k^S(\mathcal T)$ to be given by equation \eqref{eq:akS}. In Section \ref{sec:diagonal}, we prove Conjecture \ref{conj:symp-square-fund-discs} in the range where $y \ge x^2 (\log x)^{C}$ for a sufficiently large constant $C$.

We also consider the following conjecture, which is derived from \cite[Theorem 1.1]{KuperbergLalin}.
\begin{conj}\cite[Conjecture 1.9]{KuperbergLalin} \label{conj:symp-square} Let $p$ be a prime and define
\[\mathcal{S}^S_{d_k;x}(p)=\sum_{\substack{n\leq x  \\n\equiv \square \pmod{p}\\p\nmid n}}d_k(n),\]
where $\square$ indicates a perfect square.

Define the weighted variance
\begin{align*}
\mathrm{Var}^w_{p\in [y,2y]}\left(\mathcal{S}^S_{d_k;x}\right) &:= \frac{1}{y} \sum_{y < p \le 2y} \log p \Big(\mathcal{S}^S_{d_k;x}(p) - \left\langle \mathcal{S}^S_{d_k;x}\right \rangle_S\Big)^2,
\end{align*}
where $\langle \cdot \rangle_S$ denotes the limiting average; namely,
\begin{align*}
\left\langle \mathcal{S}^S_{d_k;x}\right \rangle_S&:=  \frac{1}{2y} \sum_{\substack{y < p \le 2y }} \log p \sum_{\substack{n \le x\\ p\nmid n }} d_k(n).
\end{align*}

Then, for $x^{1/k+\varepsilon}\leq y$, we have,
\begin{equation}\label{eq:conj}
\mathrm{Var}^w_{p\in [y,2y]}\left(\mathcal{S}^S_{d_k;x}\right)
 \sim a_k^{S}(\mathcal S) \frac{x}{4}\gamma_{d_k,2}^S\left(\frac{\log x}{\log y}\right) (\log y)^{2k^2+k-2},
\end{equation}
where $a_k^S(\mathcal S)$ is given by
\begin{align}\label{eq:aSkP}
a_k^S(\mathcal S)= &2\prod_p
 \left(1 -\frac{1}{p}\right)^{2k^2+k} \frac{1}{2} \left(\left(1+\frac{1}{\sqrt{p}}\right)^{-2k}+\left(1-\frac{1}{\sqrt{p}}\right)^{-2k}\right).
 \end{align}
\end{conj}
As before, the difference between \cite[Theorem 1.1]{KuperbergLalin} and Conjecture \ref{conj:symp-square} lies in the coefficient $a^S_k(\mathcal S)$. In this article, we will present a heuristic argument leading to expression \eqref{eq:aSkP}. In Section \ref{sec:symp-square}, we prove Conjecture \ref{conj:symp-square} in the range where $y \ge x^2 (\log x)^{C}$ for a sufficiently large constant $C$.

We also consider a setting studied by Rudnick and Waxman \cite{Rudnick-Waxman}. Given an ideal $\mathfrak{a}=(\alpha)$ of $\Z[i]$ one can associate a direction vector $u(\mathfrak{a})=u(\alpha):=\left(\frac{\alpha}{\overline{\alpha}}\right)^2$ in the unit circle.  Since the generator of the ideal is defined up to multiplication  by a  unit $\{\pm 1, \pm i\}$, the direction vector is well-defined. Writing
$u(\mathfrak{a})=e^{i4\theta_\mathfrak{a}}$ determines $\theta_\mathfrak{a}$, the angle of $\mathfrak{a}$, which is well-defined modulo $\pi/2$.

For a given  $\theta$, it is natural to consider the ideals $\mathfrak{a}$ such that $\theta_\mathfrak{a} \in I_K:=[\theta-\frac{\pi}{4K},\theta+\frac{\pi}{4K}]$, a small interval around $\theta$.

We have the following conjecture, which follows from \cite[Theorem 1.2]{KuperbergLalin}. 
\begin{conj}\cite[Conjecture 1.11]{KuperbergLalin2} \label{conj:symp-rudnickwaxman}
Let 
\begin{equation}\label{eq:NS-dl-definition}
\mathcal{N}^S_{d_\ell,K;x}(\theta)=\sum_{\substack{\mathfrak{a} \, \text{ideal} \\N(\mathfrak{a})\leq x\\
\theta_\mathfrak{a} \in I_K(\theta) }} d_\ell(\mathfrak{a}).
\end{equation}
Consider the variance defined by
\begin{align*}
 \mathrm{Var}\left(\mathcal{N}^S_{d_\ell,K;x}\right) := \frac 2{\pi} \int_0^{\pi/2} \big(\mathcal N^S_{d_\ell,K;x}(\theta)- \left\langle
 \mathcal{N}^S_{d_\ell,K;x} \right\rangle\big)^2 \mathrm{d}\theta,
\end{align*}
where the average is given by
\begin{equation}\label{eq:NS-dl-average}
 \left\langle \mathcal{N}^S_{d_\ell,K;x} \right\rangle:= \frac 2{\pi} \int_0^{\pi/2}  \sum_{\substack{\mathfrak{a} \, \text{ideal} \\N(\mathfrak{a})\leq x\\ \theta_\mathfrak{a} \in I_K(\theta) }} d_\ell(\mathfrak{a}) d\theta
 =\sum_{\substack{\mathfrak{a}\, \text{ideal}\\ N(\mathfrak{a})\leq x}} d_\ell(\mathfrak{a})
 \frac 2{\pi} \int_0^{\pi/2} \mathds{1}_{\theta_\mathfrak{a}\in I_K(\theta)} \mathrm{d}\theta
 = \frac{1}{K}\sum_{\substack{\mathfrak{a}\, \text{ideal}\\ N(\mathfrak{a})\leq x}} d_\ell(\mathfrak{a}).
\end{equation}

Let $x \leq K^{\ell}$. Then as $x \to \infty$, the variance $\mathrm{Var}\left(\mathcal N^S_{d_\ell,K;x}\right)$ is given by
\[\mathrm{Var}\left(\mathcal{N}^S_{d_\ell,K;x}\right)\sim a^S_\ell(\mathcal N) \frac{x}{K}
 \gamma_{d_\ell,2}^S\left(\frac{\log x}{2\log K}\right) (2\log K)^{2\ell^2+\ell-2},\]
 where
 \begin{equation}\label{eq:aSl}
a^S_\ell(\mathcal N) := 2 \left(\frac{\pi}{8b^2}\right)^{\ell(2\ell - 1)} \frac{(2+\sqrt{2})^{2\ell} + (2-\sqrt{2})^{2\ell}}{2^{\binom{2\ell + 1}{2}+1}} \prod_{p \equiv 1 \pmod 4}\left[ \left(1-\frac 1p\right)^{4\ell^2} \sum_{n=0}^\infty \binom{n+2\ell-1}{2\ell-1}^2 \frac 1{p^n}\right],
 \end{equation}
 for $b$ the Landau--Ramanujan constant given by
 \begin{equation}\label{eq:landau-ramanujan-definition}
b = \frac 1{\sqrt 2} \prod_{p \equiv 3\pmod 4} \left(1-\frac 1{p^2}\right)^{-\tfrac 12}.
 \end{equation}
 \end{conj}

In a similar way as in previous cases, the difference between \cite[Theorem 1.2]{KuperbergLalin} and Conjecture \ref{conj:symp-rudnickwaxman} lies in the coefficient $a^S_\ell(\mathcal N)$. A goal of this article is to present a heuristic argument leading to expression \eqref{eq:aSl}. In Section \ref{sec:symp-rudnickwaxman}, we prove Conjecture \ref{conj:symp-rudnickwaxman} in the range where $K \gg x (\log x)$.

 As in the unitary setting, $a^S_k(\mathcal T)$, $a^S_k(\mathcal S)$, and $a^S_\ell(\mathcal N)$ are expected to be given by Euler products. However, we notice that the expressions \eqref{eq:akS}, \eqref{eq:aSkP}, and \eqref{eq:aSl} include an extra factor of 2. As we will later discuss, this discrepancy is introduced from translating certain constraints from the number field case, which is continuous, to the discrete setting of the function field case. This feature seems to be specific for the symplectic questions that we consider in this manuscript, and it does not appear in \cite{KR3}. Curiously, this factor of $2$ does not seem to immediately appear when, for example, naively applying the ``recipe'' of \cite{CFKRS}. It is not immediately clear how to adapt the recipe to this setting. A natural formulation of the variance \eqref{eq:conj-fund-discs}, for example, is related to an average of quadratic $L$-functions via Perron's formula. As a result, the variance involves averages over quadratic $L$-functions both at the central point and all along the half-line. While the recipe allows one to relate averages of $L$-values near the central point to random matrix integrals, the relation for other points on the half line is less clear; naively extending the approximation near the central value does not produce the additional factor of $2$ in these conjectures.

 Another noticeable feature is that the constant does depend on the choice of the interval, in the sense that performing the sum over $1\leq n\leq X$ leads to a different result than performing the sum over $\alpha X\leq n \leq (\alpha+1)X$. This feature separates the number field case from the function field case, where this difference does not play a role.

 This article is organized as follows. In Section \ref{sec:leading} we obtain the leading coefficient of $\gamma_{d_k,2}^S(c)$, which we need to compare with the constants that arise from direct computation in the integer case. We recall some results due to de la Bret\`eche \cite{dlB} on the evaluation of multiple Dirichlet series in Section \ref{sec:regis-theoremes}. These results are then applied in Sections \ref{sec:diagonal}, \ref{sec:symp-square}, and \ref{sec:symp-rudnickwaxman} in order to find the respective constants $a_k^S(\mathcal{T}),a_k^S(\mathcal{S}),$ and $a_\ell^S(\mathcal{N})$ in Conjectures \ref{conj:symp-square-fund-discs}, \ref{conj:symp-square}, and \ref{conj:symp-rudnickwaxman}. Section \ref{sec:discrepancies} discusses the discrepancies due to the extra factor of 2 that is introduced in the number field case as well as the dependance on the interval for $n$. Finally, a brief discussion of numerical evidence is included in Section \ref{sec:numerics}.

 \medskip
 \noindent {\bf Acknowledgments:}  We are grateful to  Brian Conrey, Brad Rodgers, and Ezra Waxman for many helpful discussions.

\section{The leading term of $\gamma_{d_k,2}^S(c)$} \label{sec:leading}

The goal of this section is to obtain the coefficient in $\gamma_{d_k,2}^S(c)$ of the highest power of $c$ in the interval $c \in [0,1/2]$. This will allow us to compare the leading coefficient of $\gamma_{d_k,2}^S(c)$ with the constants coming from direct computation in the integer case. Recall the following result.
\begin{prop}\cite[Proposition 4.2]{KuperbergLalin2} \label{prop:nless2N}
  For $n\leq N+\frac{1+k}{2}$, we have
  \begin{equation}\label{eq:Idk2-highest-power-of-gamma}
  I_{d_k,2}^S(n;N)=\sum_{\substack{\ell=0\\\ell\equiv n \pmod{2}}}^n\binom{\frac{n-\ell}{2}+\binom{k+1}{2}-1}{\binom{k+1}{2}-1}^2\binom{\ell+k^2-1}{k^2-1}.
  \end{equation}
 Moreover, $I_{d_k,2}^S(n;N)$ is a quasi-polynomial in $n$ of period 2 and degree $2k^2+k-2$ (provided that $n\leq N+\frac{1+k}{2}$).
 \end{prop}
Our goal is to extract the coefficient of the highest power of $n$ in the right-hand side of \eqref{eq:Idk2-highest-power-of-gamma}; the highest powers of $n$  and the summation index $\ell$ in the right-hand side are equal to
\begin{align*}
&\frac{1}{(k^2-1)! \left(\left(\frac{k^2+k}{2}-1\right)!\right)^2}\sum_{\substack{\ell=0\\\ell\equiv n \pmod{2}}}^n \left(\frac{n-\ell}{2}\right)^{k^2+k-2} \ell^{k^2-1}.
\end{align*}
We proceed by dividing by $n^{2k^2 + k - 2}$ and taking the limit as $n \rightarrow \infty$ in the sum above. This gives
\[\lim_{n \rightarrow \infty} \frac{1}{n}\sum_{\substack{\ell=0\\\ell\equiv n \pmod{2}}}^n \left(\frac{1}{2}-\frac{\ell}{2n}\right)^{k^2+k-2} \left(\frac{\ell}{n}\right)^{k^2-1}=\frac{1}{2}\int_{0}^1\left(\frac{1-t}{2}\right)^{k^2+k-2} t^{k^2-1} dt,\]
 where the factor $\frac{1}{2}$ in front is due to the condition that $\ell\equiv n \pmod{2}$. In turn, the integral is given by
 \[\frac{1}{2^{k^2+k-1}}\int_{0}^1\left(1-t\right)^{k^2+k-2} t^{k^2-1} dt=\frac{(k^2+k-2)! (k^2-1)!}{2^{k^2+k-1}(2k^2+k-2)!}.\]
Putting everything together, the main coefficient of $\gamma_{d_k,2}^S(c)$ is given by
\begin{equation}
\binom{k^2+k-2}{\frac{k^2+k}{2}-1}\frac{1}{2^{k^2+k-1}(2k^2+k-2)!}.
\end{equation}

\section{Results on evaluating multiple Dirichlet series}\label{sec:regis-theoremes}

Here we record several theorems due to de la Bret\`eche \cite{dlB} on the evaluation of multiple Dirichlet series which will appear in our work later on.

Consider the general situation where $f:\mathbb{N}^\mu \rightarrow \R_{>0}$ is a positive arithmetic function. We are interested in estimating
\[S(X;\beta_1,\dots,\beta_\mu):=\sum_{1\leq n_1\leq X^{\beta_1}}\cdots \sum_{1\leq n_\mu \leq X^{\beta_\mu}} f(n_1,\dots,n_\mu).\]We will work with the Dirichlet series
\[F(s_1,\dots, s_\mu):= \sum_{n_1=1}^\infty \cdots \sum_{n_\mu=1}^\infty\frac{f(n_1,\dots,n_\mu)}{n_1^{s_1}\cdots n_\mu^{s_\mu}}.\]

Let $\mathcal{L}_\mu$ denote the set of linear forms in $(s_1,\dots,s_\mu)$ that take positive values when they are evaluated in $s_j \in \R_{>0}$. Also let $e_1,\dots,e_\mu$ denote the canonical base for $\C^\mu$.
\begin{thm}\label{thm:regis1}\cite[Th\'eor\`eme 1 (simplified version)]{dlB} Assume that the following is true:
\begin{enumerate}
 \item $F(s_1,\dots, s_\mu)$ converges absolutely in the region $\re(s_j)>\frac{1}{2}$ for $j=1,\dots, \mu$.
 \item There is a set $\{\ell^{(j)}\}_{j=1}^n$ of linear forms in $\mathcal{L}_\mu$ such that for some $\delta_1 > 0$, the function
 \[H(s_1,\dots, s_\mu):= F\left(s_1,\dots, s_\mu\right) \prod_{j=1}^n\ell^{(j)}\left(s_1-\frac{1}{2},\dots, s_\mu-\frac{1}{2}\right) \]
converges absolutely to a holomorphic function in the domain $\re(s_j)>\frac{1}{2}-\delta_1$.
\item There is a $\delta_2>0$ such that for all $\varepsilon_1, \varepsilon_2>0$, the bound
\begin{align*}|H(s_1,\dots, s_\mu)|\ll& \prod_{j=1}^n \left(\left|\im\left(\ell^{(j)}\left(s_1-\frac{1}{2},\dots, s_\mu-\frac{1}{2}\right)\right)\right|+1 \right)^{1-\delta_2\min \left\{0,\re\left(\ell^{(j)}\left(s_1-\frac{1}{2},\dots, s_\mu-\frac{1}{2}\right)\right) \right\}}\\
&\times \left(1+||\im(s_1,\dots, s_\mu)||_1^{\varepsilon_1} \right)
\end{align*}
is uniformly valid for $\re(s_j)>\frac{1}{2}-\delta_1+\varepsilon_2$.
\end{enumerate}
Then, there is a polynomial $Q_{\beta_1,\dots,\beta_\mu}(T) \in \R[T]$ of degree at most $\mu-\mathrm{rank}(\{\ell^{(j)}\}_{j=1}^n)$ and a positive real number $\theta =\theta (\mathcal{L}_h, \delta_1,\delta_2, \beta_1,\dots,\beta_\mu)$ such that for $X\geq 1$, we have
\[S(X;\beta_1,\dots,\beta_\mu)=X^{\frac{\beta_1+\cdots+\beta_\mu}{2}} (Q_{\beta_1,\dots,\beta_\mu}(\log X)+O(X^{-\theta})).\]
\end{thm}

\begin{thm}\label{thm:regis2}\cite[Th\'eor\`eme 2 (simplified version)]{dlB} Let $f$ be an arithmetic function satisfying all the hypoth\`eses from Theorem \ref{thm:regis1}. Let $\mathrm{Vect}(\{\ell^{(j)}\}_{j=1}^n)$ denote the vector space generated by the dual vectors of the linear forms in $\mathcal{L}_\mu$.
Suppose that $F(s_1,\dots,s_\mu)$ satisfies the following conditions.
\begin{enumerate}
 \item There is a function $G$ such that  $H(s_1,\dots, s_\mu)=G(\ell^{(1)}(s_1,\dots, s_\mu), \dots, \ell^{(n)}(s_1,\dots, s_\mu)))$.
 \item The vector $e_1+\cdots+e_\mu=(1,\dots,1) \in \mathrm{Vect}(\{\ell^{(j)}\}_{j=1}^n)$ and there is no proper subfamily $\mathcal{L}'$ of $\{\ell^{(j)}\}_{j=1}^n$ such that $(1,\dots,1) \in \mathrm{Vect}(\mathcal{L}')$ and
 \[\mathrm{card}(\mathcal{L}')-\mathrm{rank}(\mathcal{L}')= \mathrm{card}(\{\ell^{(j)}\}_{j=1}^n)-\mathrm{rank}(\{\ell^{(j)}\}_{j=1}^n).\]
 \end{enumerate}

Then, for $X\geq 3$, the polynomial $Q_{\beta_1,\dots,\beta_\mu}(T)$ satisfies
\[Q_{\beta_1,\dots,\beta_\mu}(\log X)=H(1/2,\dots,1/2) X^{-\frac{\beta_1+\cdots+\beta_\mu}{2}} I(X;\beta_1,\dots,\beta_\mu)+O((\log X)^{\rho-1}),\]
where $\rho=n-\mathrm{rank}(\{\ell^{(j)}\}_{j=1}^n$,
\[I(X;\beta_1,\dots,\beta_\mu):=\int_{\mathcal{A}(X;\beta_1,\dots,\beta_\mu)} \frac{dz_1\cdots dz_n}{\prod_{j=1}^n z_j^{1-\ell^{(j)}(1/2,\dots,1/2)}},\]
and
\[\mathcal{A}(X;\beta_1,\dots,\beta_\mu):= \left\{ z_1,\dots, z_n \in [1,\infty)\, : \, \prod_{j=1}^n z_j^{\ell^{(j)}(e_i)}\leq X^{\beta_i}, i=1, \dots, n\right\}.\]

\end{thm}

\begin{rem}
The original statements of Theorems \ref{thm:regis1} and \ref{thm:regis2} that can be found in \cite{dlB} allow for more freedom,
where the associated Dirichlet sum $F(s_1,\dots,s_\mu)$ is absolutely convergent in the set $\re(s_j)>\alpha_j$ for some $(\alpha_1,\dots,\alpha_\mu)\in \R_{>0}^\mu$.

We have simplified the above statements to take $(\alpha_1,\dots,\alpha_\mu)=(1/2,\dots,1/2)$, which suffices for our purposes. In addition, most of our applications will use $(\beta_1,\dots,\beta_\mu)=(1,\dots,1)$.
\end{rem}

\section{The constant in Conjecture \ref{conj:symp-square-fund-discs}}
\label{sec:diagonal}

We begin with the problem of computing the variance of the $k$-fold divisor function weighted by the quadratic character modulo $d$; namely, the problem discussed in Conjecture \ref{conj:symp-square-fund-discs}. Our goal in this section is to predict the arithmetic constant appearing in equation \eqref{eq:akS}. Our strategy is to compute the variance in a restricted range of the parameters, so that, at least conjecturally, only ``diagonal'' terms contribute.

Specifically, we begin with the variance
\begin{align}\label{eq:variance-jacobi-symbol-diagonal-expansion}
\mathrm{Var}_{r\in (y,2y]}\Bigg(\sum_{\substack{n\leq x  \\(n,r) = 1}}d_k(n)\chi_r(n)\Bigg) =& \mathbb E^*_{y < r \le 2y} \Bigg(\sum_{\substack{n\leq x  \\(n,r) = 1}}d_k(n)\chi_r(n) \Bigg)^2\\
=& \sum_{n,m \le x} d_k(n)d_k(m)\mathbb E^*_{y < r \le 2y}  \chi_r(nm), \nonumber
\end{align}
where we recall that the expectation $\mathbb E^*$ is taken over positive fundamental discriminants.

The ``diagonal'' terms are those where the expectation $\mathbb E^*_{y < r \le 2y}  \chi_r(nm)$ should be large. Since $\chi_r$ is a quadratic character, we expect this to happen precisely when $nm$ is a perfect square. Thus the diagonal terms $D$ are given by
\begin{equation*}
D = \sum_{\substack{n,m \le x \\ nm = \square}} d_k(n)d_k(m) \mathbb E^*_{y < r \le 2y}  \chi_r(nm).
\end{equation*}
If $nm$ is a square, then $\chi_r(nm)$ is $1$ if $(nm,r) = 1$, and otherwise $0$. By \cite[equation (3.1.17)]{CFKRS}, the expectation in this case is given by
\begin{equation} \label{eq:3.1.17}
\mathbb E^*_{y < r \le 2y} \chi_r(nm) \sim
\begin{cases}
\prod_{p\mid nm} \left(1 + \frac 1p\right)^{-1} & \mathrm{if}\, nm=\square,\\
o(1) & \mathrm{otherwise},
\end{cases}
\end{equation}
so that
\begin{equation*}
D \sim \sum_{\substack{n,m \le x \\ nm = \square}} d_k(n)d_k(m) \prod_{p|nm} \left(1 + \frac 1p\right)^{-1}.
\end{equation*}
In \cite{CFKRS} (3.1.17) is written as a sum over $1 \leq  r\le y$ and not $y \leq  r \le 2y$, but these are equivalent.
Note that it is crucial that the average is taken over fundamental discriminants and \emph{not} over all values $r$.

The remaining ``off-diagonal'' terms are those where the expectation should be small. If we write
\begin{equation*}
OD = \sum_{\substack{n,m \le x \\ nm \ne \square}}d_k(n)d_k(m)
\mathbb E^*_{y < r \le 2y} \chi_r(mn),
\end{equation*}
then under the generalized Riemann Hypothesis, if $y \ge x^{1/k}$,
\begin{equation*}
OD \ll x^2 y^{-1/2} \log^{C'}x
\end{equation*}
for some absolute $C'$, so for $y \ge x^2\log^{C}x$ the variance in \eqref{eq:variance-jacobi-symbol-diagonal-expansion} grows like
\begin{equation} \label{eq:var-squares}
\mathrm{Var}_{r\in (y,2y]}\Bigg(\sum_{\substack{n\leq x  \\(n,r) = 1}}d_k(n)\chi_r(n)\Bigg) \sim \sum_{\substack{n,m \le x \\ nm = \square}} d_k(n)d_k(m) \prod_{p\mid nm} \left(1 +\frac{1}{p}\right)^{-1},
\end{equation}
which we would now like to evaluate using multiple Dirichlet series techniques. Ultimately this will give us an output of the form $A_kc^{2k^2 + k - 2}$ if $\log x = c\log y$, and the constant $A_k$ in this output is precisely the arithmetic constant we would like to compute.

The diagonal terms $D$ admit a Dirichlet series
\begin{align}\label{eq:Fk}
F_k(s_1,s_2):=&\sum_{nm = \square} \frac{d_k(n)d_k(m)}{n^{s_1}m^{s_2}} \prod_{p\mid nm} \left(1+\frac 1p\right)^{-1} \nonumber\\
=& \sum_{h} \prod_{p\mid h}  \left(1+\frac 1p\right)^{-1} \sum_{n_1\dots n_k m_1\dots m_k=h^2}  \frac{1}{n_1^{s_1}\cdots n_k^{s_1} m_1^{s_2}\cdots m_k^{s_2}}.
\end{align}

Our goal is to estimate \eqref{eq:var-squares} from $F_k(s_1,s_2)$. To do this, we apply the results of de la Bret\`eche that are stated in Section \ref{sec:regis-theoremes}. Observe  from \eqref{eq:Fk} that $F_k(s_1,s_2)$ is absolutely convergent in $\re(s_1)>\frac{1}{2}$ and $\re(s_2)>\frac{1}{2}$. Moreover, writing
\begin{align*}
\mathcal{H}_k(s_1,s_2)
:=&  F_k(s_1,s_2) \zeta(2s_1)^{-\binom{k+1}{2}}\zeta(2s_2)^{-\binom{k+1}{2}}\zeta(s_1+s_2)^{-k^2},
\end{align*}
we have that
\begin{align*}
\mathcal{H}_k(s_1,s_2)=& \prod_p\left(1-\frac{1}{p^{2s_1}}\right)^{\binom{k+1}{2}}\left(1-\frac{1}{p^{2s_2}}\right)^{\binom{k+1}{2}}\left(1-\frac{1}{p^{s_1+s_2}}\right)^{k^2} \\
 &\times \left(1 + \left(1+\frac 1p\right)^{-1} \left(\tau_{s_1,\dots, s_1,s_2,\dots, s_2}(p^2)+\tau_{s_1,\dots s_1,s_2,\dots, s_2}(p^4)+\cdots \right) \right),
\end{align*}
where we define
\begin{equation}\label{eq:taudef}
\tau_{s_1,\dots, s_j}(n) := \sum_{n_1 \cdots n_j = n} \frac 1{n_1^{s_1} \cdots n_j^{s_j}}.
\end{equation}
We remark that $\mathcal{H}_k(s_1,s_2)$ converges absolutely to an analytic function in the set $\re(s_1)>\frac{1}{4}$ and $\re(s_2)>\frac{1}{4}$. In particular, for all $\epsilon > 0$, $|\mathcal H_k(s_1,s_2)| \le |\mathcal H_k(1/4 + \epsilon, 1/4+\epsilon)| \ll_{\epsilon} 1$.

One could similarly define \begin{align*}
H_k(s_1,s_2):=&  F_k(s_1,s_2) (2s_1-1)^{\binom{k+1}{2}}(2s_2-1)^{\binom{k+1}{2}}(s_1+s_2-1)^{k^2},
\end{align*}
which is the function in the statement of Theorem \ref{thm:regis1}, where $\delta_1=\frac{1}{4}$ and the family of forms is $\mathcal{L}=\{\ell^{(j_1,j_2)}(s_1,s_2)  1\leq j_1\leq j_2\leq 2k  \}$ given by
\begin{equation}\label{eq:linearforms}\ell^{(j_1,j_2)}(s_1,s_2) =\begin{cases}
                              2s_1 & 1\leq j_1\leq  j_2 \leq k,\\
                              s_1+s_2& 1 \leq j_1\leq k < j_2\leq 2k,\\
                              2s_2 & k< j_1\leq  j_2 \leq 2k.
                             \end{cases}\end{equation} Notice that
\[H_k(s_1,s_2) = \mathcal{H}_k(s_1,s_2) (2s_1-1)^{\binom{k+1}{2}}\zeta(2s_1)^{\binom{k+1}{2}} (2s_2-1)^{\binom{k+1}{2}} \zeta(2s_2)^{\binom{k+1}{2}}(s_1+s_2-1)^{k^2} \zeta(s_1+s_2)^{k^2}. \]
The bound in Theorem \ref{thm:regis1}, part (3) is then satisfied because
\begin{align*}
H_k(s_1,s_2) &= \mathcal H_k(s_1,s_2) (\zeta(2s_1)(2s_1-1))^{\binom{k+1}{2}}(\zeta(2s_2)(2s_2-1))^{\binom{k+1}{2}} (\zeta(s_1 + s_2)(s_1+s_2-1))^{k^2} \\
&\ll ((1 + |\mathrm{Im}(2s_1)|)^{1-\mathrm{Re}(2s_1)/3 + \epsilon})^{\binom{k+1}{2}}((1 + |\mathrm{Im}(2s_2)|)^{1-\mathrm{Re}(2s_2)/3 + \epsilon})^{\binom{k+1}{2}} \\
&\times ((1 + |\mathrm{Im}(s_1+s_2)|)^{1-\mathrm{Re}(s_1+s_2)/3 + \epsilon})^{k^2}.
\end{align*}

Note that
\[\tau_{1/2,\dots, 1/2}(p^{2h})=\frac{d_{2k}(p^{2h})}{p^{h}} = \binom{2h + 2k-1}{2h}\frac 1{p^h},\]
so that
\begin{align*}
 \mathcal{H}_k(1/2,1/2):=& \prod_p\left(1-\frac{1}{p}\right)^{2k^2+k} \left(1 + \left(1+\frac 1p\right)^{-1} \left( -1+\sum_{h=0}^\infty \binom{2h+2k-1}{2h}  \frac{1}{p^{h}} \right) \right)\\
   =& \prod_p\left(1-\frac{1}{p}\right)^{2k^2+k}\frac{1}{p+1} \left(1 +\frac{ p}{2}\left(\left(1+\frac{1}{\sqrt{p}}\right)^{-2k} +\left(1-\frac{1}{\sqrt{p}}\right)^{-2k}  \right) \right).
\end{align*}

Theorem \ref{thm:regis1} implies that
\begin{align*}
 \sum_{\substack{1 \le n,m \le x \\ nm = \square}}  d_k(n)d_k(m)\prod_{p\mid nm}\left(1 + \frac 1p\right)^{-1} = x(Q_{2k^2+k-2}(\log x) + O(x^{-\theta})).
\end{align*}

Our goal is to apply Theorem \ref{thm:regis2} to get an estimate of the main term. Thus, we need to verify that conditions (1) and (2) are satisfied. The variables $s_1$ and $s_2$ individually can be written as a function of the linear forms $2s_1,2s_2,$ and $s_1 + s_2$, so condition (1) holds. For condition (2), first notice that
$(1,1)\in \mathrm{Vect}(\{\ell^{(j)}\}_{j=1}^n)=\C^2$. Also, the family of forms $\mathcal L$ has $2k^2 + k$ elements and rank $2$. We need to verify that there is no proper subfamily $\mathcal{L}'$ of $\mathcal L$ such that $(1,1)\in \mathrm{Vect}(\mathcal{L}')$ and
\begin{equation}\label{eq:L}\mathrm{card}(\mathcal{L}')-\mathrm{rank}(\mathcal{L}')=2k^2+k-2.\end{equation}
The only possibility for \eqref{eq:L} to be satisfied with a proper family is for $\mathrm{rank}(\mathcal{L}')<2$. If $\mathrm{rank}(\mathcal{L}')=0$, then $\mathcal L' = \varnothing$ and \eqref{eq:L} does not hold. Thus $\mathrm{rank}(\mathcal{L}')=1$, which is only possible if all the linear forms in $\mathcal{L}'$ are of the same type. Since $(1,1)\in \mathrm{Vect}(\mathcal{L}')$, the linear forms in $\mathcal{L}'$ are all given by $s_1+s_2$. However, there are $k^2$ forms of type $s_1+s_2$, making it impossible for $\mathcal L'$ to satisfy \eqref{eq:L}.

By Theorem \ref{thm:regis2}, the polynomial $Q_{2k^2 + k-2}$ is given by
\[Q_{2k^2+k-2}(\log x) \sim \mathcal{H}_k(1/2,1/2)x^{-1} \int_{\mathcal{A}(x)} dz_{(1,1)}\cdots dz_{(2k,2k)},\]
where
\begin{align*}\mathcal{A}(x)=\left\{ z_{(j_1,j_2)}\in [1,\infty),\, 1\leq j_1\leq j_2\leq 2k\,  :\, \right. & \prod_{1\leq j_1\leq j_2 \leq k} z_{(j_1,j_2)}^2 \prod_{1 \leq j_1\leq k < j_2\leq 2k}z_{(j_1,j_2)}\leq x,\\
 & \left. \prod_{1 \leq j_1\leq k < j_2\leq 2k}z_{(j_1,j_2)}\prod_{k< j_1\leq j_2 \leq 2k} z_{(j_1,j_2)}^2\leq x\right\}.
 \end{align*}

 Our next step is evaluating the integral over $\mathcal{A}(x)$, which we accomplish by making use of the following two lemmas; both of these lemmas are proven using elementary calculus and induction on $n$.
 \begin{lem} \label{lem:calculus} Let $n \ge 1$ and $m \ge 0$ be integers. Then
\[I_{m,n}(Y):=\int_{\substack{1\leq x_1, \dots, x_n\\ x_1\cdots x_n \leq Y}} \left[\log \left( \frac{Y}{x_1\cdots x_n}\right)\right]^m \frac{dx_1\dots dx_n}{ x_1\cdots x_n}=\frac{m!}{(m+n)!} (\log Y)^{m+n}.\]
 \end{lem}
\begin{lem}
\label{lem:calculus2}
Let $n$ be a positive integer. Then
\[J_{n}(Y):=\int_{\substack{1\leq x_1, \dots, x_n\\ x_1\cdots x_n \leq Y}}dx_1\dots dx_n=(-1)^n+Y \sum_{j=0}^{n-1} \frac{(-1)^{n-1-j}}{j!} (\log Y)^j.\]
\end{lem}

Due to the particular shape of the linear forms in our problem, $I_k(X;1,1)$ is given by the volume of $\mathcal{A}(x)$.  We proceed to compute this volume. We can write
\begin{align*}
\mathrm{Vol}(\mathcal{A})= &\int_{\substack{z_{(j_1,j_2)}\in [1,\infty), 1\leq j_1\leq k <j_2\leq 2k \\ \prod_{1 \leq j_1\leq k < j_2\leq 2k}z_{(j_1,j_2)} \leq x}}   \int_{\substack{z_{(j_1,j_2)}\in [1,\infty), 1\leq j_1 \leq j_2 \leq k \\ \prod_{1 \leq j_1\leq  j_2\leq k}z_{(j_1,j_2)} \leq \left(\frac{x}{\prod_{1 \leq j_1\leq k \leq j_2\leq 2k}z_{(j_1,j_2)}}\right)^{1/2}}} \\
&\times \int_{\substack{z_{(j_1,j_2)}\in [1,\infty), k< j_1 \leq j_2 \leq 2k \\ \prod_{k <j_1\leq  j_2\leq 2k}z_{(j_1,j_2)} \leq \left(\frac{x}{\prod_{1 \leq  j_1\leq k <  j_2\leq 2k}z_{(j_1,j_2)}}\right)^{1/2}}}
dz_{(j_1,j_2)}\\
=&\int_{\substack{z_{(j_1,j_2)}\in [1,\infty), 1\leq j_1\leq k <j_2\leq 2k \\ \prod_{1 \leq j_1\leq k < j_2\leq 2k}z_{(j_1,j_2)} \leq x}}  \left( \int_{\substack{z_{(j_1,j_2)}\in [1,\infty), 1\leq j_1 \leq j_2 \leq k \\ \prod_{1 \leq j_1\leq  j_2\leq k}z_{(j_1,j_2)} \leq \left(\frac{x}{\prod_{1 \leq j_1\leq k < j_2\leq 2k}z_{(j_1,j_2)}}\right)^{1/2}}} dz_{(j_1,j_2)}\right)^2 dz_{(j_1,j_2)}.
\end{align*}
Applying Lemma \ref{lem:calculus2} to the inner integral with $n = \binom{k+1}{2}$ and $Y = \left(\frac{x}{\prod_{1 \le j_1 \le k < j_2 \le 2k} z_{(j_1,j_2)}}\right)^{1/2}$, we have
\begin{align*}
\mathrm{Vol}(\mathcal{A})
=& x \int_{\substack{z_{(j_1,j_2)}\in [1,\infty), 1\leq j_1\leq k <j_2\leq 2k \\ \prod_{1 \leq j_1\leq k < j_2\leq 2k}z_{(j_1,j_2)} \leq x}}\left(\sum_{j=0}^{\binom{k+1}{2}-1} \frac{(-1)^{\binom{k+1}{2}-1-j}}{j!2^j} \left[\log \left(\frac{x}{\prod_{1 \leq j_1\leq k < j_2\leq 2k}z_{(j_1,j_2)}}\right)\right]^j \right)^2\\&\times  \frac{dz_{(j_1,j_2)}}{\prod_{1 \leq j_1\leq k < j_2\leq 2k}z_{(j_1,j_2)}}+O\left(x^{1/2}\right)\\
=& \frac{x}{\left(\left(\binom{k+1}{2}-1\right)!\right)^22^{k^2+k-2}} \int_{\substack{z_{(j_1,j_2)}\in [1,\infty), 1\leq j_1\leq k <j_2\leq 2k \\ \prod_{1 \leq j_1\leq k < j_2\leq 2k}z_{(j_1,j_2)} \leq x}} \left[\log \left(\frac{x}{\prod_{1 \leq j_1\leq k < j_2\leq 2k}z_{(j_1,j_2)}}\right)\right]^{k^2+k-2} \\&\times  \frac{dz_{(j_1,j_2)}}{\prod_{1 \leq j_1\leq k < j_2\leq 2k}z_{(j_1,j_2)}}+O\left(x (\log x)^{2k^2+k-3}\right).
\end{align*}
For the last equality we used that the highest power of $\log x$ will come from the highest power in the sum over $j$, which can be deduced from the statement of Lemma \ref{lem:calculus}.
Applying  Lemma \ref{lem:calculus} again, we get
\begin{align*}
\mathrm{Vol}(\mathcal{A})
=&  \frac{x}{\left(\left(\binom{k+1}{2}-1\right)!\right)^22^{k^2+k-2}} \frac{(k^2+k-2)!}{(2k^2+k-2)!}\left(\log  x\right)^{2k^2+k-2}+O\left(x (\log x)^{2k^2+k-3}\right)\\
=&  \frac{x}{2^{k^2+k-2}}\binom{k^2+k-2}{\frac{k^2+k}{2}-1} \frac{1}{(2k^2+k-2)!}\left(\log  x\right)^{2k^2+k-2}+O\left(x (\log x)^{2k^2+k-3}\right).
\end{align*}

Thus we have
\begin{align*}
\sum_{\substack{1 \le n,m \le x \\ nm = \square}} d_k(n)d_k(m) \prod_{p|nm} \left(1 + \frac 1p\right)^{-1} =& \mathcal{H}_k(1/2,1/2) \frac 1{2^{k^2 + k - 2}} \binom{k^2 + k - 2}{\frac{k^2 + k}{2}-1} \frac 1{(2k^2 + k - 2)!} x (\log x)^{2k^2 + k - 2}\\& + O(x(\log x)^{2k^2 + k - 3}) \\
=& 2 \mathcal{H}_k(1/2,1/2) \gamma_k\left(\frac{\log x}{\log y}\right) x (\log y)^{2k^2 + k - 2} + O(x(\log x)^{2k^2 + k - 3}).
\end{align*}
The arithmetic constant in this case, to adjust from the ``function field'' answer of $\gamma_k\left(\frac{\log x}{\log y}\right) x (\log y)^{2k^2 + k -2}$, is $2\mathcal{H}_k(1/2,1/2)$, or
\begin{equation*}
2\prod_p\left(1-\frac{1}{p}\right)^{2k^2+k}\frac{1}{p+1} \left(1 +\frac{ p}{2}\left(\left(1+\frac{1}{\sqrt{p}}\right)^{-2k} +\left(1-\frac{1}{\sqrt{p}}\right)^{-2k}  \right) \right).
\end{equation*}
The factor of $2$ is different from the Euler product contribution, as it arises from the integral computation and not from the convergent factor $\mathcal{H}_k(s_1,s_2)$. It is also a somewhat unexpected discrepancy from the function field case; for example, no discrepancy of this form is present in the unitary case! This discrepancy is discussed in more detail in Subsection \ref{sec:factor2}.

\section{The constant in Conjecture \ref{conj:symp-square}} \label{sec:symp-square}

In this section we study a variation of the problem described in the previous section. Instead of considering the $k$-fold divisor function weighted by quadratic characters modulo $r$, we restrict $r$ to be a prime, so that the quadratic character weighting is essentially summing the $k$-fold divisor function only along quadratic residues modulo a prime $p$. Precisely, we want to understand the variance of the sum $\mathcal S^S_{d_k;x}$ defined in Conjecture \ref{conj:symp-square}. Recall that the condition $n\equiv \square \pmod{p}$ for $p\nmid n$ can be codified by $\frac{1+\chi_p(n)}{2}$ for $p\nmid n$, and therefore
\[\left\langle \mathcal{S}^S_{d_k;x}\right \rangle\sim \frac{1}{2y} \sum_{\substack{y < p \le 2y }} \log p \sum_{\substack{n \le x\\ p\nmid n }} d_k(n) .\]
 Thus, we consider the weighted variance
\begin{align*}
\mathrm{Var}_{p\in (y,2y]}^w\Bigg(\sum_{\substack{n\leq x  \\n\equiv \square \pmod{p}\\ p \nmid n}}d_k(n)\Bigg) &:= \frac 1{4y} \sum_{y < p \le 2y} \log p \Bigg(\sum_{\substack{n \le x \\ p \nmid n }}d_k(n) \chi_p(n)\Bigg)^2 \\
&= \frac 1{4y} \sum_{n,m \le x} d_k(n)d_k(m) \sum_{\substack{y < p \le 2y \\ p \nmid nm}} \log p \chi_p(nm).
\end{align*}
As in the previous section, we will ultimately only consider the setting where $y \ge x^2 \log^{C}x$, so that we may discard the off-diagonal terms. Note that in this range $p > nm$ for all terms in the sum, so the condition that $p \nmid nm$ is immaterial and we can discard it.

The diagonal terms come from the case where $nm$ is a perfect square, in which case $\chi_p(nm) = 1$ always. Thus by the prime number theorem, the diagonal terms are
\begin{align*}D:=&\frac{1}{4} \sum_{\substack{n,m\leq x\\nm=\square} } d_k(n) d_k(m) \frac 1y \sum_{\substack{y < p \le 2y}}\log p
\sim\frac{1}{4} \sum_{\substack{n,m\leq x\\nm=\square} } d_k(n) d_k(m).
\end{align*}
The off-diagonal terms are given by
\begin{align*}OD:=&\frac{1}{4} \sum_{\substack{n,m\leq x\\nm\not =\square} } d_k(n) d_k(m) \mathbb{E}_{\substack{y < p \le 2y}}\log p \chi_p(nm)
\ll  x^2 y^{-1/2}\log^{C'} x,
\end{align*}
under the generalized Riemann Hypothesis for $y \geq x^{1/k}$. Thus for $y \geq x^2 \log^{C}x$, the variance is asymptotic to the diagonal contribution.

The diagonal terms admit the Dirichlet series
\[\tilde{F}_k(s_1,s_2):=\sum_{nm=\square}\frac{d_k(n)d_k(m)}{n^{s_1}m^{s_2}}=\sum_h \sum_{n_1\cdots n_km_1\cdots m_k=h^2}\frac{1}{n_1^{s_1}\cdots n_k^{s_2}m_1^{s_1}\cdots m_k^{s_2}}.\]
Observe that $\tilde{F}_k(s_1,s_2)$ is absolutely convergent in $\re(s_1)>\frac{1}{2}$ and $\re(s_2)>\frac{1}{2}$. Writing
\[\tilde{\mathcal{H}}_k(s_1,s_2):=\tilde{F}_k(s_1,s_2)\zeta(2s_1)^{-\binom{k+1}{2}} \zeta(2s_2)^{-\binom{k+1}{2}}\zeta(s_1+s_2)^{-k^2},\]
we have that $\tilde{\mathcal H_k}(s_1,s_2)$ is very similar to the quantity $\mathcal H_k(s_1,s_2)$ in the previous section, given by
\begin{align*}
 \tilde{\mathcal{H}}_k(s_1,s_2)=& \prod_p \left(1-\frac{1}{p^{2s_1}}\right)^{\binom{k+1}{2}}\left(1-\frac{1}{p^{2s_2}}\right)^{\binom{k+1}{2}}\left(1-\frac{1}{p^{s_1+s_2}}\right)^{k^2}\\
 &\times \left(1+\tau_{s_1,\dots,s_1,s_2,\dots,s_2}(p^2)+\tau_{s_1,\dots,s_1,s_2,\dots,s_2}(p^4)+\cdots \right).
\end{align*}

We proceed as in Section \ref{sec:diagonal}, applying Theorems \ref{thm:regis1} and \ref{thm:regis2} with the linear forms family \eqref{eq:linearforms}. This eventually gives
\begin{align*}
\sum_{\substack{1 \le n,m \le x \\ nm = \square}} d_k(n)d_k(m)
=& 2 \tilde{\mathcal{H}}_k(1/2,1/2) \gamma_k\left(\frac{\log x}{\log y}\right) x (\log y)^{2k^2 + k - 2} + O(x(\log x)^{2k^2 + k - 3}),
\end{align*}
where
\begin{align*}
 \tilde{\mathcal{H}}_k(1/2,1/2)
   =& \prod_p\left(1-\frac{1}{p}\right)^{2k^2+k}\frac{1}{2}\left(\left(1+\frac{1}{\sqrt{p}}\right)^{-2k} +\left(1-\frac{1}{\sqrt{p}}\right)^{-2k} \right).
\end{align*}
As before, the factor of $2$ is different from the Euler product contribution, as it arises from the integral computation. This discrepancy is discussed in more detail in Subsection \ref{sec:factor2}.

\section{The constant in Conjecture \ref{conj:symp-rudnickwaxman}} \label{sec:symp-rudnickwaxman}

Our goal is to compute the variance
\begin{align*}
\mathrm{Var}(\mathcal N_{d_\ell,K;x}^S) &:= \frac 2{\pi} \int_0^{\pi/2} (\mathcal N_{d_\ell,K;x}^S(\theta) - \langle \mathcal N_{d_\ell,K;x}^S\rangle)^2 \mathrm d\theta
\end{align*}
in a range of $x$ and $K$ such that this computation is feasible in the integer case, in order to compare to the general conjecture in the function field case and thus predict the constant in \eqref{eq:aSl}. 
Here $\mathcal N_{d_\ell,K;x}^S(\theta)$ and $\langle \mathcal N_{d_\ell,K;x}^S\rangle$ are defined in \eqref{eq:NS-dl-definition} and \eqref{eq:NS-dl-average} respectively.

We begin by establishing the range of $x$ and $K$ that we will use in this section. Note that for two ideals $\mathfrak a = (a_1 + a_2i)$ and $\mathfrak b = (b_1 + b_2i)$ with norms bounded by $x$, where $a_1 + a_2i$ and $b_1 + b_2i$ are generators in the first quadrant so that $0 \le \mathrm{arg}(a_1 + a_2i) \le \mathrm{arg}(b_1 + b_2i) < \pi/2$, we have
\begin{align*}
\mathrm{arg}(b_1+b_2i) - \mathrm{arg}(a_1+a_2i) &= \mathrm{arg}((b_1+b_2i)(a_1-a_2i)) \\
&= \mathrm{arg}(c+di),
\end{align*}
where $c+di \in \mathbb Q[i]$ has norm $N(c+di) \le x^2$, so that $|c|\le x$ and $|d| \le x$. Moreover by the restrictions on the argument of $a_1 + a_2i$ and that of $b_1 + b_2i$, we have $0 \le \mathrm{arg}(c+di) < \pi/2$ (which also implies that $c \ne 0$). Thus
\begin{align*}
\mathrm{arg}(c+di) &= \arctan(d/c).
\end{align*}
If $d$ is $0$, then $\mathrm{arg}(c+di) = 0$; note that this occurs precisely when $\frac{a_1+a_2i}{b_1+b_2i} \in \mathbb Q$, which by slight abuse of notation we will write as $\frac{\mathfrak a}{\mathfrak b} \in \mathbb Q$. Otherwise, $d/c$ is at least $1/x$. As $x \to \infty$, $\arctan(1/x) \to 1/x$, so in particular $\arctan(1/x) \gg 1/x$. Thus $\mathrm{arg}(c+di) \gg x^{-1}$ unless $\frac{\mathfrak a}{\mathfrak b} \in \mathbb Q$.

Thus taking, say, $K \gg x \log x$, if $\mathfrak a$ and $\mathfrak b$ both appear in the sum for $\mathcal N_{d_\ell,K;x}^S(\theta)$ for the same $\theta$, then $\mathfrak a/\mathfrak b \in \mathbb Q$. These are the ``diagonal'' terms in this setting.

We now compute the variance when $K \gg x \log x$. Expanding the square, we get
\begin{align*}
\mathrm{Var}(\mathcal N_{d_\ell,K;x}^S)
= &\frac 2{\pi}\int_0^{\pi/2} \sum_{\substack{\mathfrak a, \mathfrak b \text{ ideals} \\ N(\mathfrak a),  N(\mathfrak b) \le x}} d_\ell(\mathfrak a) d_\ell (\mathfrak b) (\mathds 1_{\theta_{\mathfrak a} \in I_K(\theta)} - \tfrac 1K)(\mathds 1_{\theta_{\mathfrak b} \in I_K(\theta)} - \tfrac 1K) \mathrm d\theta.
\end{align*}
If $\mathfrak a/\mathfrak b \in \mathbb Q$, then $\theta_{\mathfrak a} = \theta_{\mathfrak b}$, so $\mathds 1_{\theta_{\mathfrak a} \in I_K(\theta)} = \mathds 1_{\theta_{\mathfrak b} \in I_K(\theta)}$. If $\mathfrak a/\mathfrak b \not\in \mathbb Q$, then $|\theta_{\mathfrak a} - \theta_{\mathfrak b}| \gg 1/x > 1/K$, so by our choice of $K$, $\mathds 1_{\theta_{\mathfrak a} \in I_K(\theta)} \mathds 1_{\theta_{\mathfrak b} \in I_K(\theta)} = 0$. Thus
\begin{align}
\mathrm{Var}(\mathcal N_{d_\ell,K;x}^S) = &\frac 2{\pi}\int_0^{\pi/2} \sum_{\substack{\mathfrak a, \mathfrak b \text{ ideals} \\ N(\mathfrak a),  N(\mathfrak b) \le x \\ \mathfrak a/\mathfrak b \in \mathbb Q }} d_\ell(\mathfrak a) d_\ell (\mathfrak b) (\mathds 1_{\theta_{\mathfrak a} \in I_K(\theta)} - \tfrac 1K)^2 \mathrm d\theta \nonumber\\
&+ \frac 2{\pi}\int_0^{\pi/2} \sum_{\substack{\mathfrak a, \mathfrak b \text{ ideals} \\ N(\mathfrak a),  N(\mathfrak b) \le x \\ \mathfrak a/\mathfrak b \not\in \mathbb Q }} d_\ell(\mathfrak a) d_\ell (\mathfrak b) (-\tfrac{1}{K}\mathds 1_{\theta_{\mathfrak a} \in I_K(\theta)}-\tfrac{1}{K}\mathds 1_{\theta_{\mathfrak b} \in I_K(\theta)} + \tfrac 1{K^2}) \mathrm d\theta \nonumber\\
= &\sum_{\substack{\mathfrak a, \mathfrak b \text{ ideals} \\ N(\mathfrak a),  N(\mathfrak b) \le x \\  \mathfrak a/\mathfrak b \in \mathbb Q }} d_\ell(\mathfrak a) d_\ell (\mathfrak b) \left( \frac 1K - \frac 1{K^2}\right) - \frac 1{K^2} \sum_{\substack{\mathfrak a, \mathfrak b \text{ ideals} \\ N(\mathfrak a),  N(\mathfrak b) \le x \\ \mathfrak a/\mathfrak b \not\in \mathbb Q }} d_\ell(\mathfrak a) d_\ell (\mathfrak b) \nonumber\\
= &\frac 1K \sum_{\substack{\mathfrak a, \mathfrak b \text{ ideals} \\ N(\mathfrak a),  N(\mathfrak b) \le x \\  \mathfrak a/\mathfrak b \in \mathbb Q }} d_\ell(\mathfrak a) d_\ell (\mathfrak b) - \frac 1{K^2} \sum_{\substack{\mathfrak a, \mathfrak b \text{ ideals} \\ N(\mathfrak a),  N(\mathfrak b) \le x }} d_\ell(\mathfrak a) d_\ell (\mathfrak b).
\label{numericssRWum}
\end{align}
We have
\begin{equation*}
\sum_{\substack{\mathfrak a \text{ ideal} \\ N(\mathfrak a) \le x}} d_\ell(\mathfrak a) \ll x(\log x)^{\ell-1},
\end{equation*}
so that the second sum is $O(x^2K^{-2}(\log x)^{2\ell-2})$. When $K \gg x \log x$, this term will be negligible compared to the size of the main term in Conjecture \ref{conj:symp-rudnickwaxman}, so we discard it. Then
\begin{align*}
\mathrm{Var}(\mathcal N_{d_\ell,K;x}^S) &\sim \frac 1K \sum_{\substack{\mathfrak a, \mathfrak b \text{ ideals} \\ N(\mathfrak a), N(\mathfrak b) \le x \\ \mathfrak a/\mathfrak b \in \mathbb Q}} d_\ell(\mathfrak a)d_\ell(\mathfrak b)= \frac 1K \sum_{\substack{\mathfrak a \text{ ideal} \\ N(\mathfrak a) \le x}} \sum_{\substack{r_1, r_2 \le \sqrt{x/N(\mathfrak a)} \\ (r_1,r_2) = 1}} d_\ell(r_1\mathfrak a)d_\ell(r_2\mathfrak a) \\
&= \frac 1K \sum_{\substack{N_1, N_2 \le x \\ N_1/(N_1,N_2) = \square \\ N_2/(N_1,N_2) = \square}} \sum_{\substack{\mathfrak a \text{ ideal} \\ N(\mathfrak a) = (N_1,N_2)}} d_\ell(\mathfrak a \sqrt{N_1/(N_1,N_2)})d_\ell(\mathfrak a \sqrt{N_2/(N_1,N_2)}).
\end{align*}
Let $\mathbb D \subseteq \mathbb Z^2$ be the subset of pairs $(N_1, N_2)$ such that $N_1/(N_1,N_2)$ and $N_2/(N_1,N_2)$ are both perfect squares.
We define $f(N_1,N_2)$ as the (complicated) inner sum, so that
\begin{equation*}
f(N_1,N_2) = \mathds{1}_{\mathbb D} \sum_{\substack{\mathfrak a \text{ ideal} \\ N(\mathfrak a) = (N_1,N_2)}} d_\ell(\mathfrak a \sqrt{N_1/(N_1,N_2)})d_\ell(\mathfrak a \sqrt{N_2/(N_1,N_2)}).
\end{equation*}
With this definition in hand, we can now consider the Dirichlet series \[F_\ell(s_1,s_2):=\sum_{N_1,N_2 \ge 1} \frac{f(N_1,N_2)}{N_1^{s_1}N_2^{s_2}}.\] This series is multiplicative, so let us look at the $p$th component for primes $p = 2$, $p \equiv 1 \pmod 4$, and $p \equiv 3 \pmod 4$.

When $p = 2$, we have (for $\mathfrak p = (1+i)$)
\begin{align*}
&\sum_{\substack{j,k = 0 \\ j \equiv k \pmod 2}}^\infty \frac {d_\ell(\mathfrak{p}^{j} )d_\ell(\mathfrak{p}^{k})}{2^{js_1 + ks_2}}=1+\frac{\ell^2}{2^{s_1+s_2}}+\frac{(\ell+1)\ell}{2}\left(\frac{1}{2^{2s_1}}+\frac{1}{2^{2s_2}}\right) +\cdots.
\end{align*}

When $p \equiv 3 \pmod 4$, we have
\begin{align*}
&\sum_{\substack{j,k = 0 \\ j,k \equiv 0 \pmod 2}}^\infty \frac {d_\ell(p^{j/2}) d_\ell(p^{k/2})}{p^{js_1 + ks_2}}= \sum_{j,k = 0}^\infty \frac{d_\ell(p^j)d_\ell(p^k)}{p^{2js_1 + 2ks_2}} = 1+\frac{\ell}{p^{2s_1}}+\frac{\ell}{p^{2s_2}}+\cdots.
\end{align*}

When $p \equiv 1 \pmod 4$, we write $(p)=\mathfrak{p}_1\mathfrak{p}_2$ for $\mathfrak p_1, \mathfrak p_2$ prime ideals in $\Z[i]$. We have
\begin{align*}
&\sum_{\substack{j,k = 0 \\ j \equiv k \pmod 2}}^\infty \frac 1{p^{js_1 + ks_2}} \sum_{\substack{\mathfrak a \subseteq \mathbb Z[i] \\ N(\mathfrak a) = p^{\min\{j,k\}} }} d_\ell( \mathfrak a p^{j-\min\{j,k\}}) d_\ell(\mathfrak a p^{k-\min\{j,k\}}) \\
&=1+\frac{2\ell^2}{p^{s_1+s_2}}+\ell^2\left(\frac{1}{p^{2s_1}}+\frac{1}{p^{2s_2}}\right)+\cdots.
\end{align*}

The Dirichlet series $F_\ell(s_1,s_2)$ is analytic for $\mathrm{Re}(s_1), \mathrm{Re}(s_2) > \frac{1}{2}$, and has poles when $\mathrm{Re}(s_1) = \frac{1}{2}$ and when $\mathrm{Re}(s_2) = \frac{1}{2}$. We consider
\begin{align*}\mathcal{H}_\ell(s_1,s_2):=&F_\ell(s_1,s_2)\zeta(2s_1)^{-\binom{\ell+1}{2}} \zeta(2s_2)^{-\binom{\ell+1}{2}}  \zeta(s_1+s_2)^{-\ell^2},
\end{align*}
which converges for $\mathrm{Re}(s_i) > \frac{1}{4}$.
We are now in a position to apply Theorems \ref{thm:regis1} and \ref{thm:regis2}. In fact, the situation is very similar to that of Section \ref{sec:diagonal}, as we also work with the family of linear forms $\mathcal{L}$ given by \eqref{eq:linearforms}. The difference in the final result lies in the evaluation of $\mathcal{H}_\ell(1/2,1/2)$. To do this, we again take advantage of the multiplicativity and proceed to consider the different components for primes $p = 2$, $p \equiv 1 \pmod 4$, and $p \equiv 3 \pmod 4$.

 For $p=2$, the contribution is given by
\begin{align*}
\left(1-\frac{1}{2}\right)^{2\ell^2+\ell}&\sum_{\substack{j,k = 0 \\ j \equiv k \pmod 2}}^\infty \frac {d_\ell(\mathfrak{p}^{j} )d_\ell(\mathfrak{p}^{k})}{2^{(j+k)/2}} = \frac 1{2^{2\ell^2+\ell}}
\sum_{\substack{j,k = 0 \\ j \equiv k \pmod 2}}^\infty
 \binom{j + \ell - 1}{\ell - 1}\binom{k + \ell - 1}{\ell - 1} \frac 1{2^{j/2}} \frac 1{2^{k/2}} \\
&= \frac 1{2^{2\ell^2+\ell}} \Bigg(\Bigg(\sum_{\substack{j = 0 \\ j \equiv 0 \pmod 2 }}^\infty \binom{j + \ell - 1}{\ell-1}\frac 1{2^{j/2}}\Bigg)^2 + \Bigg(\sum_{\substack{j = 0 \\ j \equiv 1 \pmod 2 }}^\infty \binom{j + \ell - 1}{\ell-1}\frac 1{2^{j/2}}\Bigg)^2\Bigg) \\
&= \frac{(2+\sqrt 2)^{2\ell} + (2-\sqrt 2)^{2\ell}}{2^{2\ell^2 + \ell + 1}}.
\end{align*}

For primes $p\equiv 3 \pmod{4}$, we have
\begin{align*}
\left(1-\frac{1}{p}\right)^{2\ell^2+\ell} \sum_{\substack{j,k = 0}}^\infty \frac {d_\ell(p^{j}) d_\ell(p^{k})}{p^{j+k}} =&\left(1-\frac{1}{p}\right)^{2\ell^2+\ell} \Bigg(\sum_{\substack{j = 0}}^\infty \frac {1}{p^{j}}\binom{j+\ell-1}{\ell-1}\Bigg)^2= \left(1-\frac{1}{p}\right)^{2\ell^2-\ell}.
\end{align*}

Finally, for primes $p\equiv 1\pmod{4}$, we have
\begin{align}\label{eq:p-1-mod-4-binomial-swap-all-components}
&\left(1-\frac{1}{p}\right)^{2\ell^2+\ell}\sum_{\substack{j,k = 0 \\ j \equiv k \pmod 2}}^\infty \frac 1{p^{\frac{j+k}{2}}} \sum_{r=0}^{\min\{j,k\}}d_\ell(\mathfrak p_1^r \mathfrak p_2^{\min\{j,k\}-r})d_\ell(\mathfrak p_1^{r+\frac{|k-j|}{2}} \mathfrak p_2^{\frac{|k+j|}{2}-r})\nonumber\\
&=\left(1-\frac{1}{p}\right)^{2\ell^2+\ell}\Bigg(2\sum_{j=0}^\infty \sum_{m=1}^\infty  \frac 1{p^{j+m}} \sum_{r=0}^{j}
\binom{r+\ell-1}{\ell-1}\binom{j-r+\ell-1}{\ell-1}\nonumber\\&\times \binom{r+m+\ell-1}{\ell-1}\binom{j+m-r+\ell-1}{\ell-1} + \sum_{j=0}^\infty  \frac 1{p^{j}} \sum_{r=0}^{j}
\binom{r+\ell-1}{\ell-1}^2\binom{j-r+\ell-1}{\ell-1}^2\Bigg),
\end{align}
where $2m=|k-j|$.

Consider the first sum in \eqref{eq:p-1-mod-4-binomial-swap-all-components}. We reparametrize by first replacing $m$ by $n = j + m$ and then by replacing $j$ by $t = j - r$ to get
\begin{align*}
&2 \sum_{n = 1}^\infty \frac 1{p^n}\sum_{\substack{r, t = 0 \\ r+t <  n}}^{n}\binom{r+\ell-1}{\ell-1}\binom{n-r+\ell-1}{\ell-1} \binom{t+\ell-1}{\ell-1}\binom{n-t+ \ell-1}{\ell-1}.
\end{align*}
The sum over $r,t$ remains unchanged if we make the substitutions $r \leftrightarrow n - r$ and $t \leftrightarrow n-t$, except that the limits of the sum will be that $r,t$ range from $0$ to $n$ with the constraint that $r + t > n$. If we make this change for only one copy of the sum and then combine the two copies, we get
\begin{align*}
& \sum_{n = 1}^\infty \frac 1{p^n}\sum_{\substack{r, t = 0 \\ r+t \not =  n}}^{n}\binom{r+\ell-1}{\ell-1}\binom{n-r+\ell-1}{\ell-1} \binom{t+\ell-1}{\ell-1}\binom{n-t+ \ell-1}{\ell-1}.
\end{align*}

The remaining portion of the sum in \eqref{eq:p-1-mod-4-binomial-swap-all-components} can be expressed as precisely those terms where $r + t = n$ by writing $n$ in place of $j$ and $t$ in place of $n - r$.
Thus \eqref{eq:p-1-mod-4-binomial-swap-all-components} is equal to
\begin{align*}
&\left(1-\frac 1p\right)^{2\ell^2+\ell}\sum_{n = 0}^\infty \frac 1{p^n} \left(\sum_{\substack{r= 0}}^{n}\binom{r+\ell-1}{\ell-1}\binom{n-r+\ell-1}{\ell-1}\right)^2 = \left(1-\frac 1p\right)^{2\ell^2+\ell}\sum_{n = 0}^\infty \frac 1{p^n} \binom{n + 2\ell-1}{2\ell-1}^2,
\end{align*}
where we have applied the identity
\begin{equation}
\sum_{r=0}^n \binom{r+\ell-1}{\ell-1}\binom{n-r+\ell-1}{\ell-1} = \binom{n + 2\ell-1}{2\ell-1},
\end{equation}
which is a version of the Chu--Vandermonde identity (see Equation (25) in \cite[1.2.6]{Knuth}) and it holds for all $r, n,$ and $\ell$.

Putting all of this together, we finally get
\begin{align*}
\mathcal{H}_\ell(1/2,1/2)&= \frac{(2+\sqrt 2)^{2\ell} + (2-\sqrt 2)^{2\ell}}{2^{2\ell^2 + \ell + 1}} \prod_{p\equiv 3 \pmod{4}} \left(1-\frac{1}{p}\right)^{2\ell^2-\ell}\\ &\times  \prod_{p\equiv 1 \pmod{4}} \Bigg[\left(1-\frac{1}{p}\right)^{2\ell^2+\ell} \sum_{n=0}^\infty \binom{n+2\ell -1}{2\ell-1}^2\frac{1}{p^n}\Bigg].
\end{align*}
The above expression can be simplified further, since
\begin{align*}
 \prod_{p\equiv 3 \pmod{4}} \left(1-\frac{1}{p}\right)^{2\ell^2-\ell} \prod_{p\equiv 1 \pmod{4}}\left(1-\frac{1}{p}\right)^{2\ell^2+\ell}
 &=\left(\frac{\pi}{8b^2}\right)^{2\ell^2-\ell}\prod_{p\equiv 1 \pmod{4}}\left(1-\frac{1}{p}\right)^{4\ell^2},
\end{align*}
where  $b$ is the Landau--Ramanujan constant defined in \eqref{eq:landau-ramanujan-definition} and we have made use of the Euler product expression for $L(1,\chi_{-4}) = \pi/4$. Thus we finally obtain 
\begin{align*}
\frac 1K \sum_{\substack{\mathfrak a, \mathfrak b \text{ ideals} \\ N(\mathfrak a), N(\mathfrak b) \le x \\ \mathfrak a/\mathfrak b \in \mathbb Q}} d_\ell(\mathfrak a)d_\ell(\mathfrak b)= 2 \mathcal{H}_\ell(1/2,1/2) \gamma_\ell\left(\frac{\log x}{2\log K}\right) \frac xK (2\log K)^{2\ell^2 + \ell - 2}+ O(x(\log x)^{2\ell^2 + \ell - 3}),
\end{align*}
where
\begin{align*}
\mathcal{H}_\ell(1/2,1/2)=&\left(\frac{\pi}{8b^2}\right)^{\ell(2\ell-1)} \frac{(2+\sqrt{2})^{2\ell}+(2-\sqrt{2})^{2\ell}}{2^{2\ell^2+\ell+1}} \prod_{p\equiv 1 \pmod{4}}\Bigg[\left(1-\frac{1}{p}\right)^{4\ell^2}\sum_{n=0}^\infty \binom{n+2\ell-1}{2\ell -1}^2 \frac{1}{p^n}\Bigg].
\end{align*}
Once again, the factor of $2$ is different from the rest of the constant, as it arises from the integral computation. This discrepancy is discussed in more detail in the following section.

\section{Discrepancies from the function field case and factors not present in the unitary case} \label{sec:discrepancies}

We seem to be noticing two types of discrepancy between $\gamma_{d_k,2}^S(c)$ and the heuristic answers in the integer case (beyond the arithmetic constant given by an Euler product).

\subsection{A factor of $2$} \label{sec:factor2}
Here we will discuss the factor of $2$ discrepancy appearing in all of the problems we study. This effect is already visible in the case when $k = 1$, so for simplicity we will restrict to this case. We will also only discuss the variance of the $k$-fold divisor function weighted by $\chi_r(n)$, although the same factor is present for the same reasons in the other cases.

When $k = 1$, many of the sums we get in the function field case and in the integer case are roughly the same. In particular, in the computation using de la Bret\`eche's estimates of multiple Dirichlet series given in Theorems \ref{thm:regis1} and \ref{thm:regis2}, the function $\gamma_1(c)$ is given by $X^{-1}\mathrm{Vol}(\mathcal{A}(X,X,X))$ where
\begin{equation*}
\mathcal{A}(X,X,X) = \{1 \le y_1,y_2,y_3 < \infty : y_1 y_2^2 \le X, y_1 y_3^2 \le X\}.
\end{equation*}
As computed in Section \ref{sec:diagonal}, the volume of $\mathcal A(X,X,X)$ is $\sim X \log X$.

Considering the analogous volume in the function field case, however, we get a smaller answer by a factor of $\tfrac 12$. Here we would consider
\begin{equation*}
\mathcal A_{\mathbb F_q}(N,N,N) = \{f_1,f_2,f_3 \in \mathbb F_q[T] \text{ monic} : |f_1f_2^2| = q^N, |f_1f_3^2| = q^N\},
\end{equation*}
where for $f \in \mathbb F_q[T]$, $|f|:=q^{\deg(f)}$.
Note that the constraint that $y_1y_2^2 \le X$ has become the constraint that $|f_1f_2^2| = q^N$; in the large $q$ limit, almost all polynomials of degree $\le N$ have degree exactly $N$, so this is equivalent.

However, in this case, we have
\begin{align*}
\#\mathcal A_{\mathbb F_q}(N,N,N) &= \sum_{\substack{d_1,d_2,d_3 \le N \\ d_1 + 2d_2 = N \\ d_1 + 2d_3 = N}} q^{d_1+d_2+d_3} = \sum_{\substack{d_1,d_2 \le N \\ d_1+2d_2 = N}} q^N = \left\lceil\frac N2\right\rceil q^N. 
\end{align*}
The parity condition on $d_1$ is the source of the final factor of $\frac 12$ in the computation of $\gamma_{d_k,2}^S(c)$, which does not appear in the integer case.

Notably, this discrepancy is not present in problems with unitary symmetry types.

\subsection{The choice of interval}

Another type of discrepancy comes from the fact that the answer in the integer case depends on the interval chosen for the variance, not just the size of the interval. To see what this means in practice, let us consider the problem of the variance of $\sum_{Ax < n \le Bx} d_k(n) \chi_r(n)$, as $r$ varies over fundamental discriminants. Two natural versions of this problem to consider are when $A = 0$ and $B = 1$, or when $A = 1$ and $B = 2$. Notably, there is heuristic and numeric evidence that the answer in the integer case depends on the choice of $A$ and $B$ (rather than depending only on $B-A$, for example).

We begin by expanding
\begin{equation*}
\mathbb{E}^*_{y<r\leq 2y}\Bigg(\sum_{\substack{Ax < n \le Bx\\(n,r)=1}} d_k(n)\chi_r(n)\Bigg)^2 = D + ND,
\end{equation*}
where $D$ denotes diagonal terms and $ND$ denotes off-diagonal terms.

The diagonal terms are given by
\begin{equation*}
D = \sum_{\substack{Ax < n,m \le Bx \\ nm = \square}} d_k(n)d_k(m) \mathbb{E}^*_{y<r\leq 2y} \chi_r(nm).
\end{equation*}
If $nm$ is a square, then by a result of Jutila (see for example, \cite[equation (3.1.21)]{CFKRS}),
\begin{equation*}
D \sim \sum_{\substack{Ax < n,m \le Bx \\ nm = \square}} d_k(n)d_k(m) \prod_{p|nm} \left(1 + \frac 1p\right)^{-1}.
\end{equation*}

Meanwhile the off-diagonal terms are
\begin{equation*}
OD = \sum_{\substack{Ax < n,m \le Bx \\ nm \ne \square}}d_k(n)d_k(m)
\mathbb{E}^*_{y<r\leq 2y} \chi_r(nm).
\end{equation*}
Under GRH, if $y \ge x^{1/k}$,
\begin{equation*}
OD \ll x^2 y^{-1/2} \log^Cx
\end{equation*}
for some absolute $C$, so for $y \ge x^2\log^{C'}x$ the variance grows like
\begin{equation*}
\mathrm{Var}_{r\in (y,2y]}\Bigg(\sum_{\substack{Ax < n\leq Bx  \\(n,r) = 1}}d_k(n)\chi_r(n)\Bigg)  \sim \sum_{\substack{Ax < n,m \le Bx \\ nm = \square}} d_k(n)d_k(m) \prod_{p|nm} \left(1 + \frac1p\right)^{-1},
\end{equation*}
which we can evaluate using multiple Dirichlet series techniques and Theorems \ref{thm:regis1} and \ref{thm:regis2} as in Section \ref{sec:diagonal}. This gives
\begin{align*}
 \sum_{\substack{Ax < n,m \le Bx \\ nm = \square}}  d_k(n)d_k(m)\prod_{p\mid nm}\left(1 + \frac 1p\right)^{-1} = x(Q_{2k^2+k-2}(\log x) + O(x^{-\theta})),
\end{align*}
and Theorem \ref{thm:regis2} then tells us that
\[Q_{2k^2+k-2}(\log x) \sim \mathcal{H}_k(1/2,1/2)x^{-1} \int_{\mathcal{A}_{[A,B]}(x)} dz_{(1,1)}\cdots dz_{(2k,2k)},\]
where
\begin{align*}\mathcal{A}_{[A,B]}(x)=\left\{ z_{(j_1,j_2)}\in [1,\infty),\, 1\leq j_1\leq j_2\leq 2k\,  :\, \right. & Ax \le \prod_{1\leq j_1\leq j_2 \leq k} z_{(j_1,j_2)}^2 \prod_{1 \leq j_1\leq k < j_2\leq 2k}z_{(j_1,j_2)}\leq Bx,\\
 & \left. Ax \le \prod_{1 \leq j_1\leq k < j_2\leq 2k}z_{(j_1,j_2)}\prod_{k< j_1\leq j_2 \leq 2k} z_{(j_1,j_2)}^2\leq Bx\right\}.
 \end{align*}
The volume of $\mathcal{A}_{[A,B]}(x)$ will depend on $A$ and $B$. To compute it, we proceed as in Section \ref{sec:diagonal}, getting that
\begin{align*}
\mathrm{Vol}(\mathcal{A}_{[A,B]})
=&\int_{\substack{z_{(j_1,j_2)}\in [1,\infty), 1\leq j_1\leq k <j_2\leq 2k \\ \prod_{1 \leq j_1\leq k < j_2\leq 2k}z_{(j_1,j_2)} \in [Ax,Bx]}} \\& \times  \Big( \int_{\substack{z_{(j_1,j_2)}\in [1,\infty), 1\leq j_1 \leq j_2 \leq k \\ \prod_{1 \leq j_1\leq  j_2\leq k}z_{(j_1,j_2)} \in \left[\left(\frac{Ax}{\prod_{1 \leq j_1\leq k \leq j_2\leq 2k}z_{(j_1,j_2)}}\right)^{1/2}, \left(\frac{Bx}{\prod_{1 \leq j_1\leq k \leq j_2\leq 2k}z_{(j_1,j_2)}}\right)^{1/2}\right]}} dz_{(j_1,j_2)}\Big)^2 dz_{(j_1,j_2)}.
\end{align*}
Applying Lemma \ref{lem:calculus2} to the inner integral with $n = \binom{k+1}{2}$ and $Y = \left(\frac{Bx}{\prod_{1 \le j_1 \le k < j_2 \le 2k} z_{(j_1,j_2)}}\right)^{1/2}$ and substracting the case where $B$ is replaced by $A$, we have
\begin{align*}
\mathrm{Vol}(\mathcal{A}_{[A,B]})
=& x \int_{\substack{z_{(j_1,j_2)}\in [1,\infty), 1\leq j_1\leq k <j_2\leq 2k \\ \prod_{1 \leq j_1\leq k < j_2\leq 2k}z_{(j_1,j_2)} \in [Ax,Bx]}}\left(\sum_{j=0}^{\binom{k+1}{2}-1} \frac{(-1)^{\binom{k+1}{2}-1-j}}{j!2^j} \left\{\left[B^{1/2}\log \left(\frac{Bx}{\prod_{1 \leq j_1\leq k < j_2\leq 2k}z_{(j_1,j_2)}}\right)\right]^j\right.\right.\\ &\left. \left.- \left[A^{1/2}\log \left(\frac{Ax}{\prod_{1 \leq j_1\leq k < j_2\leq 2k}z_{(j_1,j_2)}}\right)\right]^j\right\}  \right)^2  \frac{dz_{(j_1,j_2)}}{\prod_{1 \leq j_1\leq k < j_2\leq 2k}z_{(j_1,j_2)}}+O\left(x^{1/2}\right)\\
=& \frac{(\sqrt B -\sqrt A)^2x}{\left(\left(\binom{k+1}{2}-1\right)!\right)^22^{k^2+k-2}} \int_{\substack{z_{(j_1,j_2)}\in [1,\infty), 1\leq j_1\leq k <j_2\leq 2k \\ \prod_{1 \leq j_1\leq k < j_2\leq 2k}z_{(j_1,j_2)} \in  [Ax,Bx]}} \left[ \log \left(\frac{x}{\prod_{1 \leq j_1\leq k < j_2\leq 2k}z_{(j_1,j_2)}}\right)\right]^{k^2+k-2}\\& \times   \frac{dz_{(j_1,j_2)}}{\prod_{1 \leq j_1\leq k < j_2\leq 2k}z_{(j_1,j_2)}}+O\left(x (\log x)^{2k^2+k-3}\right).
\end{align*}
For the last step, we took the highest power of $\log x$ in the sum over $j$, which  has the largest contribution by the statement of Lemma \ref{lem:calculus}. Finally, applying  Lemma \ref{lem:calculus} again, and keeping the terms with the highest power of $\log x$, we get
\begin{align*}
\mathrm{Vol}(\mathcal{A}_{[A,B]})
=&  \frac{(\sqrt B-\sqrt A)^2x}{\left(\left(\binom{k+1}{2}-1\right)!\right)^22^{k^2+k-2}} \frac{(k^2+k-2)!}{(2k^2+k-2)!}\left(\log  x\right)^{2k^2+k-2}+O\left(x (\log x)^{2k^2+k-3}\right)\\
=&  \frac{(\sqrt B-\sqrt A)^2x}{2^{k^2+k-2}}\binom{k^2+k-2}{\frac{k^2+k}{2}-1} \frac{1}{(2k^2+k-2)!}\left(\log  x\right)^{2k^2+k-2}+O\left(x (\log x)^{2k^2+k-3}\right).
\end{align*}

Thus if $A=0$ and $B = 1$, we get a main term of $2 x\gamma_{d_k,2}^S\left(\log x\right)$, whereas if $A =1$ and $B = 2$, we get a main term of $(\sqrt 2 - 1)^2 2 x\gamma_{d_k,2}^S\left(\log x\right)$, and if $A = 1/2$ and $B = 3/2$, we get a main term of $(\sqrt 3 - 1)^2 x\gamma_{d_k,2}^S\left(\log x\right)$. This discrepancy cannot be present in the function field case in the limit as $q \to \infty$, since as $q \to \infty$ almost all polynomials of degree at most $x$ have degree exactly $x$.

\section{Numerics}\label{sec:numerics}

We present the results of some computations that we have done in support of our conjectures. Each graph depicts, for the fixed value of $c$ and one of our conjectures, the ratio between the numerical computation and our expected result for different values of $x$, for $k=1$ and for $k=2$. Note that the scales for each graph are different, as they depend on technical limitations. These numerics are severely restricted by the growth rate of the computation time involved, so that $x$ is not sufficiently large in these computations to appreciate the convergence towards the expected value 1. While the graphs for $k=1$ suggest that our conjectures are correct, the computations for $k=2$ merely suggest that our conjectures recover the correct rate of growth, but are not sufficient to draw conclusions about the arithmetic constants involved.

\begin{figure}
\begin{tabular}{cc}
\includegraphics[width=18pc]{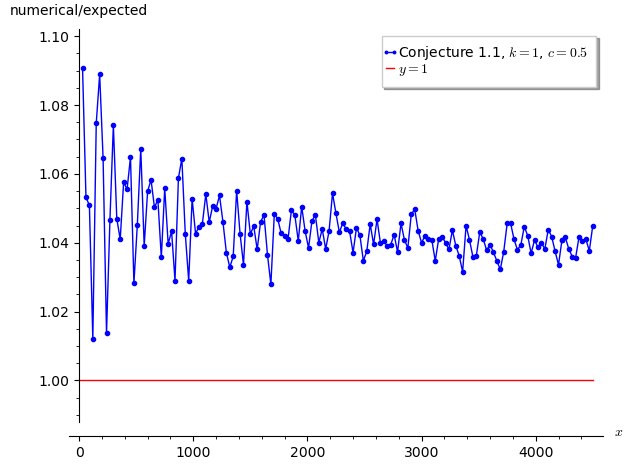}&
\includegraphics[width=18pc]{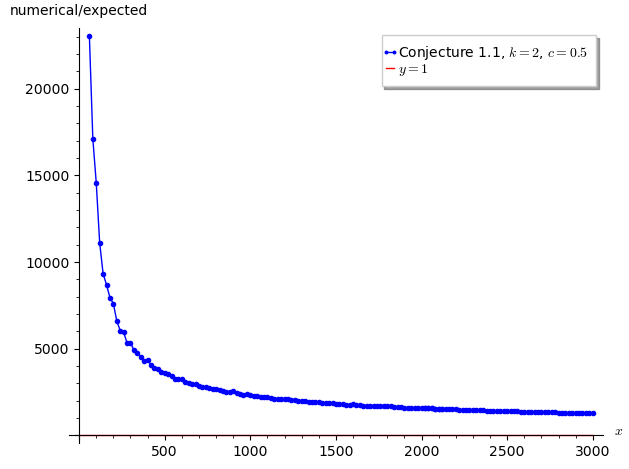}
\end{tabular}
\caption{\label{fig:conj1.1} Graph of ratios between the numerical output and our expectation for the variance in Conjecture \ref{conj:symp-square-fund-discs} when $\frac{\log x}{\log y} = 0.5$ for the cases $k=1$ (left) and $k=2$ (right).}
\end{figure}

\begin{figure}
\begin{tabular}{cc}
\includegraphics[width=18pc]{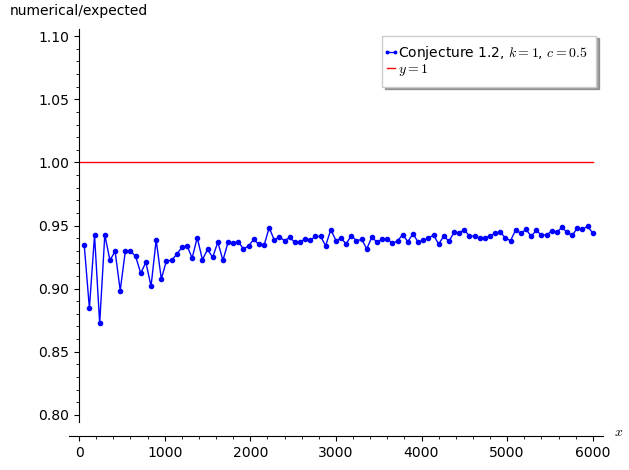}&
\includegraphics[width=18pc]{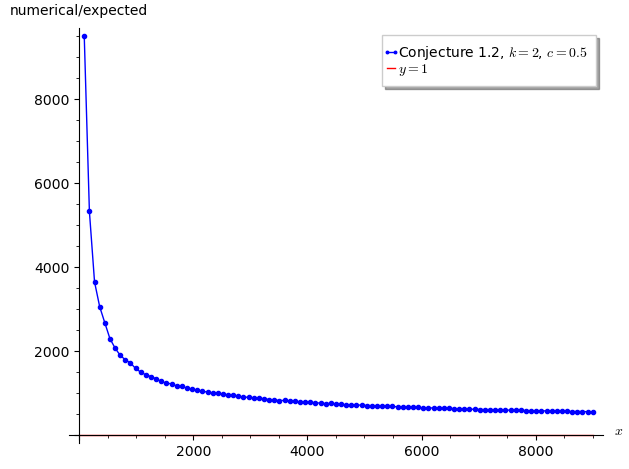}
\end{tabular}
\caption{\label{fig:conj1.2}Graph of ratios between the numerical output and our expectation for the variance in Conjecture \ref{conj:symp-square} when $\frac{\log x}{\log y} = 0.5$ for the cases $k=1$ (left) and $k=2$ (right).}
\end{figure}

Figures \ref{fig:conj1.1} and \ref{fig:conj1.2} present Conjectures \ref{conj:symp-square-fund-discs} and  \ref{conj:symp-square} respectively for $c=\frac{\log x}{\log y}=0.5$. When $k=1$, the ratios between the numerical computation and our expected result are quite close to 1 (with a rough error of 5\%). When $k=2$, the ratio is disminishing fast and appears to approach an asymptotic, but the value of this asymptotic is ambiguous.


\begin{figure}
\begin{tabular}{cc}
\includegraphics[width=18pc]{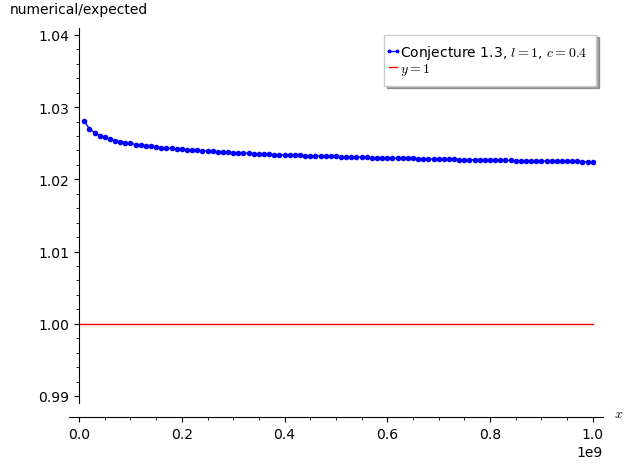}&
\includegraphics[width=18pc]{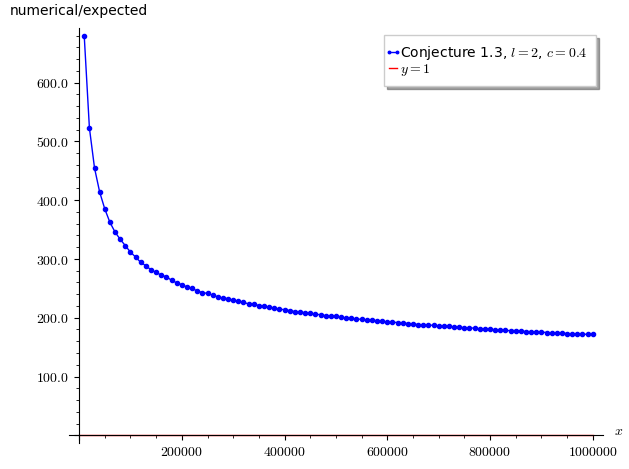}
\end{tabular}
\caption{\label{fig:conj1.3} Graph of ratios between the numerical output and our expectation for the variance in Conjecture \ref{conj:symp-rudnickwaxman} when $\frac{\log x}{2\log K} = 0.4$ for the cases $\ell=1$ (left) and $\ell=2$ (right).}
\end{figure}

\begin{figure}
\begin{tabular}{cc}
\includegraphics[width=18pc]{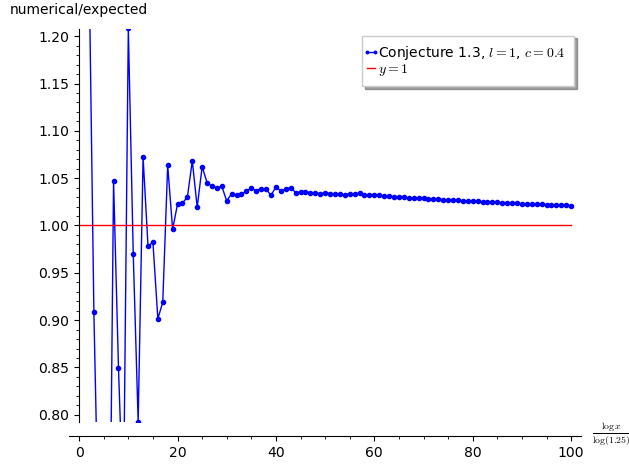}&
\includegraphics[width=18pc]{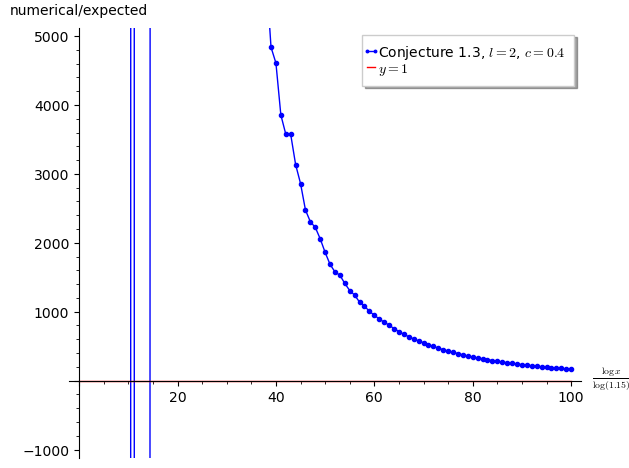}
\end{tabular}
\caption{\label{fig:conj1.3-logarithmic}  Graph of ratios between the numerical output and our expectation for the variance in Conjecture \ref{conj:symp-rudnickwaxman} when $\frac{\log x}{2\log K} = 0.4$ with logarithmic scale
for the cases $\ell=1$, basis $1.25$  (left) and $\ell=2$, basis $1.15$ (right).}
\end{figure}

Figure \ref{fig:conj1.3} depicts the ratio corresponding to Conjecture
 \ref{conj:symp-rudnickwaxman} when $c = \frac{\log x}{2\log K} = 0.4$. In this range (but not when $c = 0.5$), the variance is given by equation \eqref{numericssRWum}.  When $\ell = 1$, the values are again quite close to $1$; these values appear to approach an asymptotic slightly higher than $1$ on the linear scale, but in Figure \ref{fig:conj1.3-logarithmic}, they are depicted on a logarithmic scale, and seem to clearly continue descending towards $1$. Figure \ref{fig:conj1.3-logarithmic} depicts a noisy regime when $x$ is small and a reasonably smooth approximation for larger $x$; in this case this corresponds to the fact that the ``diagonal terms'' in \eqref{numericssRWum} is precisely the variance when $K \gg x \log x$, which happens when $c = 0.4$ only for large enough $x$. Meanwhile the plot for $\ell = 2$ resembles the other conjectures for $k=2$ in which the graph seems to approach an asymptotic with unclear constant.

\bibliographystyle{amsalpha}

\bibliography{Bibliography}

\end{document}